\documentstyle[12pt]{article}
\def\Ran{\mathop{\rm Ran}\nolimits}

\def\Re{\mathop{\rm Re}\nolimits}

\def\diag{\mathop{\rm diag}\nolimits}

\def\e{{\varepsilon}}

\def\b{{\beta}}

\def\a{{\alpha}}
\def\d{{\delta}}

\def\l{{\lambda}}

\def\D{{\Delta}}

\def\o{{\omega}}
\def\q{{\qquad\rule{.4em}{.5em}}}
\def\be{\begin{equation}}
\def\ee{\end{equation}}
\newtheorem{thm}{Theorem}
\newtheorem{lem}[thm]{Lemma}
\newtheorem{prop}[thm]{Proposition}
\newtheorem{dfn}[thm]{Definition}

\newcommand{\norm}[1]{\left\| #1 \right\|}

\date{February 5, 1997}
\author{V.~M.~Manuilov}
\title{An invariant for pairs of almost commuting unbounded operators}
\begin{document}

\maketitle

 \begin{abstract}
For a wide class of pairs of unbounded selfadjoint operators $(A,B)$ with
bounded commutator we construct a $K$-theoretical integer invariant
$\o(A,B)$, which is continuous, is equal to zero for commuting operators
and is equal to one for the pair $(x,-i\,d/dx)$.

 \end{abstract}

For a wide class of pairs of unbounded selfadjoint operators $(A,B)$ with
bounded commutator and with operator $(1+A^2+B^2)^{-1}$ being compact we
construct a $K$-theoretical integer invariant $\o(A,B)$,
which is continuous, is equal to zero for commuting operators and is equal
to one for the pair $(x,-i\,d/dx)$.
Such invariants in the case of {\em bounded} operators were constructed
by T.~Loring \cite{e-l,lor} for matrices realizing `almost commutative'
torus and sphere.

Let $A$ and $B$ be two unbounded selfadjoint operators acting on a
separable Hilbert space $H$ with dense common domain $D(A)\cap D(B)$.
We suppose also that their compositions $AB$ and $BA$
have dense common domain and such that their commutator $AB-BA$ is
bounded on $H$.
Remark that as $D(A)\cap D(B)$ is dense in $H$, so the operators $A\pm iB$
are densely defined and are formally adjoint to each other, hence these
operators are closable. Denote by $C$ the closure of the operator
$A+iB$, $C=(A-iB)^*$.
Put $\D=1+(A+iB)(A-iB)$. $\D$ is a selfadjoint invertible operator and
we suppose that its inverse operator $\D^{-1}$ is {\em compact}.

Remember that the group $K^0_{comp}({\bf R}^2)=K_0(C_0({\bf R}^2))$ is
generated by the Bott element which can be defined by a projection-valued
function of the variable $z=x+iy$ on ${\bf C}$:
 \begin{equation}\label{Bott}
P=P(x,y)=\left(
\begin{array}{cc}
{\displaystyle\frac{1}{1+|z|^2}} &
{\displaystyle\frac{z}{1+|z|^2}} \smallskip \\
{\displaystyle\frac{\overline{z}}{1+|z|^2}} &
{\displaystyle 1-\frac{1}{1+|z|^2}}
\end{array}
\right);
 \end{equation}
$P^2=P$.
One can look at this function as at the matrix of an operator acting on
the Hilbert space $H\oplus H$ (made out of the operators of multiplication
by $x$ and by $y$). Change now the operators $x$ and $y$ in (\ref{Bott})
by the operators $A$ and $B$:
 $$
Q=Q(A,B)=\left(
\begin{array}{cc}
\D^{-1} & C\D^{-1}\\
\D^{-1}C^* & 1-\D^{-1}
\end{array}
\right).
 $$
As $A$ and $B$ do not commute, there are different ways to substitute $A$
and $B$ instead of $x$ and $y$; we have fixed one of them.
We will show that if
$\norm{[A,B]}$ is small enough then the operator $Q$ is close to
a projection, namely for any $\e>0$ there exists some $\d>0$ such that if
$\norm{[A,B]}<\d$ then $\norm{Q^2-Q}<\e$.

 \begin{lem}
The operators $C(\l+C^*C)^{-1}$ and $(\l+CC^*)^{-1}C$ are bounded for any
$\l>0$ and $\norm{C(\l+C^*C)^{-1}}\leq \frac{1}{\sqrt{\l}}$,
$\norm{(\l+CC^*)^{-1}C}\leq \frac{1}{\sqrt{\l}}$.

 \end{lem}

{\bf Proof.} The statement follows from the estimates
 $$
\norm{C(\l+C^*C)^{-1/2}}^2=
\norm{(\l+C^*C)^{-1/2}C^*C(\l+C^*C)^{-1/2}}=\sup_{t\geq
0}\frac{t}{\l+t}=1,
 $$
 $$
\norm{C(\l+C^*C)^{-1}}\leq
\norm{C(\l+C^*C)^{-1/2}}\norm{(\l+C^*C)^{-1/2}}\leq
\frac{1}{\sqrt{\l}}.
 $$
Taking $C^*$ instead of $C$ we obtain the estimate
 \begin{equation}\label{soprjazh}
\norm{C^*(\l+CC^*)^{-1/2}}\leq \frac{1}{\sqrt{\l}}
 \end{equation}
and taking an adjoint in (\ref{soprjazh}) we get
 $$
\norm{(\l+CC^*)^{-1/2}C}\leq \frac{1}{\sqrt{\l}}.\q
 $$

 \begin{lem}\label{commut} One has
$C(1+C^*C)^{-1}=(1+CC^*)^{-1}C$.

 \end{lem}

{\bf Proof.}
As the operator $(1+C^*C)^{-1}$ is compact, so the Hilbert
space $H$ contains the dense linear subspace $\cup_{k=1}^\infty H_k$
where $H_k={\rm Span}\,\{\xi_1,\ldots,\xi_k\}$, $\xi_k$ being the
(ordered) eigenvectors of $(1+C^*C)^{-1}$.

Let $\{\l_k\}$ be the ordered (increasing) set
of eigenvalues of the operator $C^*C$. Let $\xi\in H_k$ for some
$k$ and suppose that $\l\in{\bf C}$, $\Re\l>\l_k$. Then we have
 \begin{eqnarray}\label{series}
C(\l+C^*C)^{-1}\xi&=&\l^{-1}C\left(1+\frac{1}{\l}C^*C\right)^{-1}\xi
\nonumber\\
&=&\l^{-1}C\left(\xi-\frac{1}{\l}C^*C\xi+\frac{1}{\l^2}(C^*C)^2\xi-
\ldots\right)
\\
&=&\l^{-1}\left(\xi-\frac{1}{\l}CC^*C\xi+\frac{1}{\l^2}(CC^*)^2C\xi-
\ldots\right).\nonumber
 \end{eqnarray}
By definition
 \begin{equation}\label{domain}
D(C^*C)=\{\psi\in D(C): C\psi\in D(C^*)\}.
 \end{equation}
As the operator $C$ is closed, so the operator $C^*C$ is selfadjoint
and $\Ran C^*C=H$ (see
\cite{reed-simon}), therefore the eigenvectors
$\xi_k$ lie in the domain of the operator $C^*C$.
Hence $\xi\in D(C^*C)$. But it follows from (\ref{domain}) that
$\xi\in D(C)$, $C\xi\in D(C^*)$. So the vector $C\xi$ lies in the domain
of operators of the form $(CC^*)^n$.
As the space $H_k$ is
finite-dimensional, so the operator $C$ is bounded on $H_k$,
$\norm{C|_{H_k}}<c_k$ for some $c_k$.
The norm of the $n$-th term in (\ref{series}) is not greater than
$c\norm{\xi}\l_k^n/|\l|^n$
and the series (\ref{series}) is convergent when $\Re\l>\l_k$.
So we have
 \begin{eqnarray*}
C(\l+C^*C)^{-1}\xi&=&
\l^{-1}\left(1-\frac{1}{\l}CC^*+\frac{1}{\l^2}(CC^*)^2-\ldots\right)C\xi\\
&=&(\l+CC^*)^{-1}C\xi.\nonumber
 \end{eqnarray*}

Consider two analytic operator-valued functions
 $$
f_1(\l)=C(\l+C^*C)^{-1}, \quad
f_2(\l)=(\l+CC^*)^{-1}C,\quad \Re\l>0.
 $$
They are bounded and they coincide on $H_k$ for $\Re\l$ big
enough, hence they coinside on $H_k$ on the whole domain of the functions
$f_i(\l)$, i.e. for $\Re\l>0$. But as $\cup_{k=1}^\infty H_k\subset H$ is
dense, so $f_1(\l)=f_2(\l)$ on $H$. \q

 \begin{lem}\label{estim-commut}
If $\norm{C^*C-CC^*}<\e$ then one has
 \begin{eqnarray*}
&&\norm{\left((1+C^*C)^{-1}-(1+CC^*)^{-1}\right)C}\leq
\frac{\e}{2(1+\e)},\\
&&\norm{C\left((1+C^*C)^{-1}-(1+CC^*)^{-1}\right)}\leq
\frac{\e}{2(1+\e)}.
 \end{eqnarray*}

 \end{lem}

{\bf Proof}\, of the first assertion follows from the estimate
 \begin{eqnarray*}
\norm{\left((1{+}C^*C)^{{-}1}{-}(1{+}CC^*)^{{-}1}\right)C}
&=&\norm{\left((1{+}C^*C)^{{-}1}{-}(1{+}C^*C{+}X)^{{-}1}\right)C}\\
&=&\norm{\left((1{+}C^*C)^{{-}1}{-}
(1{+}(1{+}C^*C)^{{-}1}X)^{{-}1}(1{+}C^*C)^{{-}1}\right)C}\\
&=&\norm{\left(1{-}(1{+}(1{+}C^*C)^{{-}1}X)^{{-}1}\right)(1{+}C^*C)^{{-}1}C}\\
&\leq&\left(1-\frac{1}{1-\e}\right)\norm{(1+C^*C)^{-1}C}\leq\frac{\e}{(1-\e)},
 \end{eqnarray*}
where $X=CC^*-C^*C$.
The second assertion can be proved in the same way.  \q

 \begin{lem}\label{neprer}
Let $E,F$ be selfadjoint positive operators such that
the operator $E-F$ is densely defined on $H$ and
$\norm{E-F}<\e$.
Put ${\displaystyle f(t)=\frac{t}{(1+t)^2}}$. Then one has
 $$
\norm{f(E)-f(F)}<\frac{3\e-\e^2}{(1-\e)^2}.
 $$

 \end{lem}

{\bf Proof.} Put $E-F=X$, $(1+F)^{-1}X=Y$, $(1+Y)^{-1}-1=Z$.
One has $\norm{Y}<\e$, ${\displaystyle \norm{Z}<\frac{\e}{1-\e}}$.
Then
 \begin{eqnarray*}
E(1{+}E)^{{-}2}{-}F(1{+}F)^{{-}2}&{=}
&(F{+}X)(1{+}F{+}X)^{{-}2}{-}F(1{+}F)^{{-}2}\\
&{=}&(F{+}X)\left((1{+}Y)^{{-}1}(1{+}F)^{{-}1}
(1{+}Y)^{{-}1}(1{+}F)^{{-}1}\right){-}F(1{+}F)^{{-}2}\\
&{=}&(F{+}X)\left((1{+}Z)(1{+}F)^{{-}1}(1{+}Z)
(1{+}F)^{{-}1}\right){-}F(1{+}F)^{{-}2}\\
&{=}&F\left(Z(1{+}F)^{{-}1}{+}(1{+}F)^{{-}1}Z{+}
Z(1{+}F)^{{-}1}Z(1{+}F)^{{-}1}\right)\\
&{+}&X\left((1{+}F)^{{-}2}{+}Z(1{+}F)^{{-}1}
{+}(1{+}F)^{{-}1}Z{+}Z(1{+}F)^{{-}1}Z(1{+}F)^{{-}1}\right)\\
&{\leq}&2\frac{\e}{1-\e}+\left(\frac{\e}{1-\e}\right)^2+\frac{\e}{(1-\e)^2}
=\frac{3\e-\e^2}{(1-\e)^2}.\q
 \end{eqnarray*}

Now we are ready to show that for small enough $\e$ the operator $Q$ is
close to a projection.
 \begin{thm}
If $\norm{C^*C-CC^*}<\e$ then
${\displaystyle\norm{Q^2-Q}<\frac{4\e-2\e^2}{(1-\e)^2}}$.

 \end{thm}
{\bf Proof.}
Easy calculations (using lemma \ref{commut}) show that \quad $Q^2-Q=$
 $$
\left(
\begin{array}{cc}
(1{+}C^*C)^{{-}2}C^*C{-}(1{+}CC^*)^{{-}2}CC^*&
C((1{+}CC^*)^{{-}1}{-}(1{+}C^*C)^{{-}1})(1{+}C^*C)^{{-}1}\\
(1{+}C^*C)^{{-}1}((1{+}CC^*)^{{-}1}{-}(1{+}C^*C)^{{-}1})C^*&
0
\end{array}
\right).
 $$
But the norm of the off-diagonal elements is estimated by the lemma
\ref{estim-commut}:
 \begin{eqnarray*}
&&\norm{C((1+CC^*)^{-1}-(1+C^*C)^{-1})(1+C^*C)^{-1}}\\
&&\leq
\norm{C((1+CC^*)^{-1}-(1+C^*C)^{-1})}\norm{(1+C^*C)^{-1}}
\leq
\frac{\e}{1-\e}
 \end{eqnarray*}
and the norm of the first element is estimated by the lemma \ref{neprer}:
 $$
\norm{(1+C^*C)^{-2}C^*C-(1+CC^*)^{-2}CC^*}\leq \frac{3\e-\e^2}{(1-\e)^2}.
 $$
From the inequality
 $$
\norm{\left(
\begin{array}{cc}
\a&\b\\
\b^*&0
\end{array}
\right)}\leq \norm{\a}+\norm{\b}
 $$
one finally has
 $$
\norm{Q^2-Q}<\frac{3\e-\e^2}{(1-\e)^2}+\frac{\e}{1-\e}=
\frac{4\e-2\e^2}{(1-\e)^2}=\e'.\q
 $$

So we see that the point $1/2$ divides the spectrum of $Q(A,B)$ into two
subsets when $\norm{[A,B]}<0.02$.

Denote by $H_N\subset H$ the subspace generated by the vectors
$e_{N+1},\ldots$ of some basis $\{e_i\}$ of $H$, $H=L_N\oplus H_N$.
Our `almost projection' $Q$ acts on the Hilbert space $H\oplus H$.
Consider the decomposition of this space:
 $$
H\oplus H=(L_N\oplus L_N)\oplus (H_N\oplus H_N).
 $$
With respect to this decomposition the operator $Q$ has the form
 $$
Q=
\left(
\begin{array}{cc}
q_{11}&q_{12}\\
q_{12}^*&q_{22}
\end{array}
\right).
 $$
It follows now from compactness of $\D^{-1}$ that for big enough $N$
$\norm{q_{12}}$ is close to zero and $q_{22}$ is close to the operator
 \begin{equation}\label{proj01}
\left(
\begin{array}{cc}
0 & 0\\
0 & 1
\end{array}
\right).
 \end{equation}
Therefore the operator $q_{11}$ should be close to a projection,
$\norm{q_{11}^2-q_{11}}<\e''$ for some $\e''$ and the eigenvalues of
$q_{11}$ also satisfy $|\l^2-\l|<\e''$.
When $\e''<1/4$ then all
eigenvalues of the operator $q_{11}$ can be divided into two sets:
 $$
S_0=\{\l\in\ {\rm Spec}\,q_{11}:\l<1/2\}\quad {\rm and}\quad
S_1=\{\l\in\ {\rm Spec}\,q_{11}:\l>1/2\}.
 $$
 \begin{dfn}{\rm
Let $M_N=\#(S_1)$ be the number of
eigenvalues close to one
of the operator $q_{11}$ when $\e',\e''<1/4$. Denote by $\o(A,B)$ the
number $M_N-N$.
}
 \end{dfn}
 \begin{prop}
The definition of the number $\o(A,B)$ depends neither on the choice of a
basis in $H$ nor on the number $N$.

 \end{prop}
{\bf Proof.}
Notice that when $\e',\e''<1/4$ then the spectrum of restriction of the
operator $Q$ to the subspace $L_N\oplus L_N$ does not contain the point
$1/2$.
For big enough $N,N'$, $N<N'$ the restriction of $Q$ to the
space $(L_{N'}\ominus L_N)\oplus(L_{N'}\ominus L_N)$ is close to the
projection of the form
(\ref{proj01}), so on $(L_{N'}\ominus L_N)\oplus(L_{N'}\ominus L_N)$ the
number of
eigenvalues less than $1/2$ is equal to the number of eigenvalued bigger
than $1/2$, therefore $\o(A,B)$ does not depend on $N$. If we take two
bases $\{e_i\}$ and $\{e'_i\}$ then the space $L'_N$ generated by
$e'_1,\ldots,e'_N$ can be approximated by the space $L_{N'}$ for big
enough $N'$.\q

 \begin{thm}
Let $A$ and $B$ commute. Then $\o(A,B)=0$.

 \end{thm}

{\bf Proof.} Remember that by definition (see \cite{reed-simon})
$A$ and $B$ commute when their spectral projections commute.
It follows from the spectral theorem for two commuting operators and from
the compactness of $\D^{-1}$ that it would be enough to consider the case
of operators acting on $H=L^2({\bf Z}\times{\bf Z})$ by multiplication:
 $$
A(f)_{nm}=a_nf_{nm},\quad B(f)_{nm}=b_mf_{nm},
 $$
where $f\in H$. Then
the matrix of the projection $Q$ is $2{\times}2$-block-diagonal of the
form:
 $$
Q=\diag\{q_{nm}\} \quad {\rm where} \quad q_{nm}=
\left(
\begin{array}{cc}
{\displaystyle\frac{1}{a_n^2+b_m^2+1}} &
{\displaystyle\frac{a_n+ib_m}{a_n^2+b_m^2+1}} \medskip \\
{\displaystyle\frac{a_n-ib_m}{a_n^2+b_m^2+1}} &
{\displaystyle 1-\frac{1}{a_n^2+b_m^2+1}}
\end{array}
\right).
 $$
But each $2{\times}2$-matrix $q_{nm}$ is a projection unitarily
equivalent to
$\left(
\begin{array}{cc}
0&0\\0&1
\end{array}
\right)$,
hence the numbers of eigenvalues equal to zero and equal to one coincide
on $L_N\oplus L_N$ for every $N$.\q

Show now that the topological invariant $\o(A,B)=0$ is continuous.
Let $\{A_n\}$ be a sequence of closed symmetric operators such that
for any $\phi\in H$ one has
 $$
\norm{(A_n-A)\phi}\leq a_n\norm{A\phi}+b_n\norm{\phi}.
 $$
We say that this sequence tends to the operator $A$ if the sequences
$\{a_n\}$ and $\{b_n\}$ tend to zero. By the Kato -- Rellich theorem
\cite{reed-simon} the operators $A_n$ are selfadjoint on the domain of $A$
when $|a_n|<1$.

 \begin{thm}\label{continuity}
Let $A_n\to A$, $B_m\to B$ and let the sequence of commutators $[A_n,B_m]$
converge to $[A,B]$. If $\norm{[A,B]}<0.02$
then $\o(A_n,B_m)=\o(A,B)$ for big enough $n$ and $m$.

 \end{thm}
{\bf Proof.} We start with proving continuity with respect to the first
argument. The second one can be treated in the same way. Denote
$A_n-A=X_n$ and fix $n$ so that $a_n<\d$, $b_n<\d$, $[X_n,B]<\d$.
Then
 \begin{equation}\label{normX}
\norm{X_n\phi}<\d(\norm{A\phi}+\norm{\phi}).
 \end{equation}
It follows from (\ref{normX}) that
 $$
\norm{X_n(1+A^2)^{-1/2}\phi}<\d
\left(\norm{A(1+A^2)^{-1/2}\phi}+\norm{(1+A^2)^{-1/2}\phi}\right)
\leq 2\d,
 $$
hence $\norm{X_n(1+A^2)^{-1/2}}<2\d$.
As $1+A^2\leq 1+(A+iB)(A-iB)+0.02=\D+0.02$, so
$\D^{-1}\leq\frac{1}{1-0.02}(1+A^2)^{-1}<2(1+A^2)^{-1}$ and
 $$
\norm{\D^{-1/2}X_n}^2=\norm{X_n\D^{-1}X_n}
\leq 2\norm{X_n(1+A^2)^{-1}X_n}=2\norm{(1+A^2)^{-1/2}X_n}^2,
 $$
therefore
 $$
\norm{\D^{-1/2}X_n}<2\d; \quad
\norm{X_n\D^{-1/2}}<2\d.
 $$
Then
 $$
\norm{\D^{-1/2}AX_n\D^{-1/2}}<\d,\quad
\norm{\D^{-1/2}X_nA\D^{-1/2}}<\d,\quad
\norm{\D^{-1/2}X_n^2\D^{-1/2}}<4\d^2.
 $$
Estimate the matrix elements of the operator $Q(A_n,B)-Q(A,B)$.
For diagonal elements we have
 \begin{eqnarray*}
&&\norm{(1+(A+X_n+iB)(A+X_n-iB))^{-1}-(1+(A+iB)(A-iB))^{-1}}\\
&&=\norm{\D^{{-}1/2}}\cdot
\norm{\left(1+\D^{-1/2}
(AX_n+X_nA+X^2+i[X_n,B])\D^{-1/2}
\right)^{-1}-1}
\cdot\norm{\D^{-1/2}}
\\
&&
<\frac{1}{1-(2\d+4\d^2+\d)}-1;
 \end{eqnarray*}
which implies also compactness of the operator
$(1+(A+X_n+iB)(A+X_n-iB))^{-1}$.
A similar estimate holds for the off-diagonal elements:
 \begin{eqnarray*}
&&\norm{(A{+}X_n{+}iB)(1{+}(A{+}X_n{+}iB)(A{+}X_n{-}iB))^{-1}-
(A{+}iB)(1{+}(A{+}iB)(A{-}iB))^{-1}}\\
&&{\leq}\norm{X_n\D^{-1/2}}
\cdot
\norm{
\left(1+\D^{-1/2}(AX_n+X_nA+X_n^2+i[X_n,B])
\D^{-1/2}\right)^{-1}}
\cdot\norm{
\D^{-1/2}}\\
&&{+}\norm{(A{+}iB)\D^{-1/2}}
{\cdot}
\norm{
\left(
\left(1{+}\D^{-1/2}(AX_n{+}X_nA{+}X_n^2{+}i[X_n,B])
\D^{-1/2}\right)^{-1}
{-}1\right)}
{\cdot}\norm{
\D^{-1/2}}
\\
&&
<2\d\frac{1}{1-(3\d+4\d^2)}+\left(\frac{1}{1-(3\d+4\d^2)}-1\right).
 \end{eqnarray*}
So the norms of the matrix elements of $Q(A_n,B)-Q(A,B)$ tend to zero,
hence the sequence of operators $Q(A_n,B)-Q(A,B)$ tends to zero too.
\q

Remark that our definition of convergence implies that for $\a_0\neq 0$
one has $\a A\to \a_0 A$. Notice also that compactness of the operator
$(1+(A+iB)(A-iB))^{-1}$ implies compactness of the operators $(1+(\l
A+i\mu B)(\l A-i\mu B))^{-1}$ for all $\l>0$, $\mu>0$.
It makes the following definition correct.

 \begin{dfn}\label{d-om}
{\rm
Let $A$, $B$ be selfadjoint operators such that
\begin{enumerate}
\item
$AB$ and $BA$ have dense common domain in $H$ and $[A,B]$ is
bounded on $H$,
\item
$A$ and $B$ have dense common domain in $H$ and
$(1+(A+iB)(A-iB))^{-1}$ is compact.
\end{enumerate}
Then put $\o(A,B)=\o(\l A,\mu B)$ for any $\l>0$ and $\mu>0$
such that $\norm{[\l A,\mu B]}$ is small enough.
}
 \end{dfn}

Notice that instead of projection (\ref{Bott}) we could take any
projection realizing some class of $K^0_{comp}({\bf R}^2)$. Then a pair
of operators $(A,B)$ gives a map $K^0_{comp}({\bf R})\rightarrow {\bf Z}$,
hence pairs of operators represent $K$-homology classes of ${\bf R}^2$.

As a corollary of the theorem \ref{continuity} we see that if $T$ is a
bounded selfadjoint operator on $H$ commuting with $B$,
then $\o(A+T,B)=\o(A,B)$ for any pair
$(A,B)$. In particular, if $A$ and $B$ commute then $\o(A+T,B)=0$.
Remark that if the spectrum of $A$ or $B$ has a lacuna of big enough
length then by method of \cite{manUMN4,manFA} one can also obtain
$\o(A,B)=0$. It reflects the topological analogy: the reduced
$K$-homology group of ${\bf R}^2$ with a strip cutted out is trivial.
Formulas
 $$
h_1=(1+\D)^{-1},\quad h_2=A(1+\D)^{-1},\quad h_3=B(1+\D)^{-1}
 $$
can be viewed as the formulas of the stereographic projection. They define
a `non-commutative sphere' given by three compact selfadjoint operators
satisfying relations $\norm{h_1^2+h_2^2+h_3^2-h_1}<\e$ for some $\e>0$
\cite{lor}.

 \begin{thm}
Let $A=x$, $B=-i\frac{d}{dx}$. Then $\o(A,B)=1$.

 \end{thm}
{\bf Proof.} Multiply the operators $A$ and $B$ by $\sqrt{\l}$ to make
their commutator small enough. As $\norm{[A,B]}=2\l$, so we can take $\l$
so that ${\displaystyle\frac{16\l-32\l^2}{(1-4\l)^2}<1/4}$
(then $\norm{Q^2-Q}<1/4$). It
will be easier to make calculations in the basis of $H$ consisting of
eigenvectors $\phi_n$ of the operator
$C^*C=\l(x^2-d^2/dx^2-1)$.
Then one has
 $$
C\phi_n=\sqrt{n\l}\phi_{n-1},\quad
C^*\phi_n=\sqrt{(n+1)\l}\phi_{n+1},\quad
\D\phi_n=(2n\l+1)\phi_n,
 $$
and the operator $Q$ can be written in the form
 $$
\left(
\begin{array}{ccccccccc}
1&0&0&0&0&0&0&0&\cdots\\
0&0&\frac{\sqrt{\l}}{2\l+1}&0&0&0&0&0&\cdots\\
0&\frac{\sqrt{\l}}{2\l+1}&\frac{1}{2\l+1}&0&0&0&0&0&\cdots\\
0&0&0&\frac{2\l}{2\l+1}&\frac{\sqrt{2\l}}{4\l+1}&0&0&0&\cdots
\vspace{1pt}\\
0&0&0&\frac{\sqrt{2\l}}{4\l+1}&\frac{1}{4\l+1}&0&0&0&\cdots\\
0&0&0&0&0&\frac{4\l}{4\l+1}&\frac{\sqrt{3\l}}{6\l+1}&0&\cdots
\vspace{1pt}\\
0&0&0&0&0&\frac{\sqrt{3\l}}{6\l+1}&\frac{1}{6\l+1}&0&\cdots\\
0&0&0&0&0&0&0&\frac{6\l}{6\l+1}&\cdots\\
\cdots&\cdots&\cdots&\cdots&\cdots&\cdots&\cdots&\cdots&\cdots
\end{array}
\right).
 $$
To calculate the number of eigenvalues close to one we notice first that
the matrix $q_{11}$ (which is the restriction of $Q$ to the space
$L_N\oplus L_N$) decomposes into direct sum of $N-1$
two-dimensional blocks of the form
 \begin{equation}\label{2x2}
\left(
\begin{array}{cc}
\frac{2(n-1)\l}{2(n-1)\l+1}&\frac{\sqrt{n\l}}{2n\l+1}\smallskip\\
\frac{\sqrt{n\l}}{2n\l+1}&\frac{1}{2n\l+1}
\end{array}
\right)
 \end{equation}
and of two one-dimensional blocks: $1$ and $\frac{2N\l}{2N\l+1}$. When $N$
is big enough the last eigenvalue is close to one:
$\lim_{N\to\infty}\frac{2N\l}{2N\l+1}=1$. So it remains to show that out
of the two eigenvectors of each matrix of the form (\ref{2x2}) one is
close to zero and another is close to one.
Let $x^2-px+q$ be the characteristic polynomial for the matrix
(\ref{2x2}). Then
 $$
1-p=\frac{2\l}{(2n\l+1)(2(n-1)\l+1)}\leq\frac{2\l}{2\l+1}<2\l,
 $$
 $$
q=\frac{2\l}{(2n\l+1)^2(2(n-1)\l+1)}\leq\frac{2\l}{(2\l+1)^2}<2\l,
 $$
hence when $\l$ tends to zero the characteristic polynomial tends to
$x^2-x$ and that gives us necessary estimates for the eigenvalues, so
every block of the form (\ref{2x2}) gives equal numbers of eigenvalues in
the sets $S_0$ and $S_1$. But the two scalar blocks give us $M_N=N+1$.\q

\bigskip
Remark that the construction of the $K$-theoretical invariant can be
generalized to the case of pairs of operators on Hilbert $C^*$-modules
over $C^*$-algebras. This generalization includes the case of families
of pairs of operators $(A_x,B_x)$, $x\in X$ satisfying pointwise the
conditions of definition \ref{d-om}. We refer to \cite{lance} for
basic facts about Hilbert $C^*$-modules.
Let $H_{\cal A}$ be a Hilbert ${\cal A}$-module and let $(A,B)$ be a
pair of unbounded selfadjoint operators on $H_{\cal A}$ such that
\begin{enumerate}
\item
compositions $AB$ and $BA$ have dense common domain in $H_{\cal A}$ and
the operator $AB-BA$ is bounded,
\item
$A$ and $B$ have dense common domain in $H_{\cal A}$ and
the operator $(k+(A+iB)(A-iB))^{-1}$ is ${\cal A}$-compact on $H_{\cal A}$.
\end{enumerate}
Then the invariant $\o(A,B)$ can be defined as an element of
$K_0({\cal A})$ given by $\o(A,B)=M_N-N$ where $M_N$ is the $K$-theory
class of the spectral projection of $Q$ corresponding to the interval
$(-\infty,1/2)$.

\bigskip
{\bf Acknowledgement.}  The present paper was prepared with the partial
support of RBRF (grant N 96-01-00276). I am grateful to M.~Frank,
A.~S.~Mishchenko and E.~V.~Troitsky for helpful discussions.


\vspace{.8cm}
\noindent
V.~M.~Manuilov\\
Dept. of Mech. and Math.,\\
Moscow State University,\\
Moscow, 119899, RUSSIA\\
e-mail: manuilov@mech.math.msu.su

\end{document}